\documentclass[12pt]{article}
\usepackage{amsfonts,amssymb,amsbsy,amsmath,amsthm,enumerate}

\usepackage{epsfig}
\usepackage{color}

\numberwithin{equation}{section}

\newtheorem{theorem}{Theorem}[section]
\newtheorem{lemma}{Lemma}[section]
\newtheorem{follow}{Corollary}[section]
\newtheorem{pr}{Proposition}[section]
\theoremstyle{definition}

\newcommand{\Subsection}[1]{\subsection{ #1} ${}^{}$}

\newcommand{\D}{{\mathcal D}}

\newcommand{\Z}{{\mathbb Z}}

\newcommand{\R}{{\mathbb R}}

\newcommand{\C}{{\mathbb C}}
\newcommand{\N}{{\mathbb N}}

\newcommand{\res}{\operatorname{Res}}

\newcommand{\im}{\operatorname{Im}}

\def\phi {\varphi}

\newcommand{\cH}{{\mathcal H}}
\newcommand{\cZ}{{\mathcal Z}}
\newcommand{\cN}{{\mathcal N}}

\newcommand{\xp}{X_\perp}
\newcommand{\mper}{m_\perp}

\def\<{\langle}
\def\>{\rangle}

\newcommand{\bel}{\begin{equation} \label}
\newcommand{\ee}{\end{equation}}

\newcommand{\hopa}{H_{\parallel}}

\newcommand{\hope}{H_{\perp}}

\newcommand{\rt}{{\mathbb R}^{3}}
\newcommand{\rd}{{\mathbb R}^{2}}
\newcommand{\re}{{\mathbb R}}
\newcommand{\ce}{{\mathbb C}}

\newcommand{\fract}[2]{\genfrac{}{}{0pt}{}{\scriptstyle #1}{\scriptstyle #2}}

\begin{document}
\begin{center}{\LARGE \bf  Resonances and Spectral Shift Function Singularities for Magnetic
Quantum Hamiltonians}

\medskip

{\sc Jean-Fran\c{c}ois Bony, Vincent Bruneau, Georgi Raikov}

\medskip
\end{center}

\bigskip

{\bf Abstract.} In this survey article we consider the operator pair $(H,H_0)$
where $H_0$ is the shifted 3D Schr\"odinger operator with
constant magnetic field, $H : = H_0 + V$, and $V$ is  a short-range electric potential  of a fixed sign. We describe the asymptotic behavior  of the
Krein spectral shift function (SSF) $\xi(E; H,H_0)$ as the energy $E$ approaches the Landau levels $2bq$, $q \in {\mathbb Z}_+$, which play the role of
thresholds in the spectrum of $H_0$. The main asymptotic term of $\xi(E; H,H_0)$ as $E \to 2bq$ with a fixed $q \in {\mathbb Z}_+$ is written in terms of appropriate compact Berezin--Toeplitz operators. Further, we investigate the relation between the threshold singularities of the SSF and the accumulation of resonances at the Landau levels. We establish the existence of resonance--free sectors adjoining any given Landau level and prove that the number of the resonances in the complementary sectors is infinite. Finally, we obtain the main asymptotic term of the local resonance counting function near an arbitrary fixed Landau level; this main asymptotic term is again expressed via the Berezin--Toeplitz operators which govern the asymptotics of the SSF at the Landau levels.

\bigskip

\bigskip

\bigskip

\section{Introduction}
\label{s1}
This article is a survey of our results on two closely related topics:
\begin{itemize}
\item The threshold singularities of the Krein spectral shift function (SSF) for the operator pair $(H,H_0)$ where $H_0$ is the shifted 3D Schr\"odinger operator with constant magnetic field, and $H: = H_0 + V$ with an appropriate short-range potential $V$ of a fixed, positive or negative, sign;
    \item The accumulation of the resonances of the operator $H$ at its spectral thresholds.
    \end{itemize}
 The self-adjoint unperturbed operator
    \bel{1}
H_0 = H_0(b) : = (-i\nabla - A)^2 - b
    \ee
is defined initially on $C_0^{\infty}(\re^3)$ and then is closed in $L^2(\re^3)$. Here
   $A =
\left(-\frac{bx_2}{2}, \frac{bx_1}{2}, 0\right)$ is a
magnetic potential generating the magnetic field
    $B = {\rm curl}\;A = (0,0,b)$ where $b>0$ is the constant scalar intensity of $B$. \\
It is well known that the spectrum $\sigma(H_0)$ of the operator $H_0$ is absolutely continuous
 and coincides with $[0,\infty)$ (see e.g. \cite{ahs}). Moreover, the so called {\em Landau levels} $2bq$, $q
\in {\mathbb Z}_+ : = \{0,1,2,\ldots\}$, play the role of  thresholds in
$\sigma(H_0)$ (see below Subsection \ref{ss22a}).\\
The perturbation of $H_0$ is the electric potential $V: {\rt} \to \re$. Assume that $V$ is  Lebesgue measurable, and bounded. On the domain of $H_0$ define the perturbed operator
    \bel{2}
H : = H_0 + V
    \ee
 which obviously is self-adjoint in $L^2(\rt)$.\\  In order to describe the assumptions on the decay of $V$ we need the following notations. For $x \in \rt$ we will occasionally write $x = (\xp, x_3)$ where $\xp = (x_1,x_2) \in \rd$ are the variables in the plane perpendicular to the magnetic field, and $x_3 \in \re$ is the variable along the magnetic field.
We will suppose that $V$ satisfies one of the following estimates:
\begin{itemize}
\item D (anisotropic decay): $V(x) = O
(\langle\xp\rangle^{-m_{\perp}}\langle x_3 \rangle^{-m_3})$ with $m_{\perp} > 2$,
$m_3 > 1$;

\item ${\rm D}_0$ (isotropic decay): $V(x) = O(\langle x \rangle^{-m_0})$ with $m_0>3$;

\item ${\rm D}_{\rm exp}$ (fast decay with respect to $x_3$): $V(x) = O
(\langle\xp\rangle^{-m_{\perp}}\exp{(-N|x_3|)})$ with some $m_{\perp} >
0$ and any $N > 0$.
\end{itemize}
Note that assumption ${\rm D}_0$ implies D. Moreover, evidently, assumption ${\rm D}_{\rm exp}$  with $m_{\perp} > 2$ again implies D.\\
The article is organized as follows. In Section \ref{s2} we discuss
the behavior of the SSF for the operator pair $(H , H_0)$ near the
Landau levels $2bq$, $q \in \Z_+$. Its main asymptotic term (see
Theorem \ref{th1}) is given in terms of auxiliary Berezin--Toeplitz
operators discussed in detail in Subsection \ref{ss22}. As an example
of the possible applications of Theorem \ref{th1}, a generalized
Levinson formula is deduced from it (see Corollary \ref{f1}). In
Subsection \ref{ss27} we describe briefly the extensions of Theorem \ref{th1} to Pauli and Dirac operators with non constant magnetic field. \\
 Further, Section \ref{s3} is dedicated to  the resonances  and the embedded eigenvalues of $H$. In Subsection \ref{ss31} we state Theorem \ref{thsan1} which implies that under very general assumptions about the decay of $V$ the operator $H$ has infinitely many eigenvalues embedded in its essential spectrum, provided that $V$ is axisymmetric and non positive. Next, in Subsection \ref{ss32}  we define the resonances of $H$ as the poles of a meromorphic continuation of the resolvent $(H-z)^{-1}$, ${\rm Im}\,z > 0$, to an appropriate infinitely sheeted Riemann surface. Further, we establish the existence of resonance-free regions and of regions with infinitely many resonances in a vicinity of each Landau levels (see Theorem \ref{th2}). Finally, for every $q \in \Z_+$ fixed, we obtain in Theorem \ref{th3} the main asymptotic term of the number of the resonances on an annulus centered at the Landau level $2bq$ as its inner radius tends to zero. \\
 We have been working on the problems discussed in the article for almost a decade. Some of the results obtained have already been surveyed (see e.g \cite{bbbr, r2}). We decided to include them again in the present opus since we wanted to tell here our story from the very beginning,  referring the reader when necessary to the original works but preferably not to other surveys.

\section{Singularities of the spectral shift function at the Landau levels}
\label{s2}

\subsection{The spectral shift function $\xi(E;H,H_0)$}
\label{ss21}

Let $V$ satisfy D. Then the diamagnetic inequality easily implies that the operator $|V|^{1/2}(H_0+1)^{-1}$  is Hilbert--Schmidt, and hence the resolvent difference $(H-i)^{-1} - (H_0-i)^{-1}$ is a
trace-class operator.  Therefore, there exists a unique
$$ \xi =
 \xi(\cdot;H,H_0)  \in L^1(\re; (1+E^2)^{-1}dE)
$$
 such that
the Lifshits--Krein trace formula
$$
{\rm Tr}\,(f(H)-f(H_0)) = \int_{\re} \xi(E; H,H_0) f'(E) dE
$$
holds for each $f \in C^{\infty}_0(\re)$ and the normalization condition
$\xi(E;H,H_0) = 0$ is fulfilled for each
 $E \in (-\infty,\inf \sigma(H))$ (see the original works \cite{lif, kr} or \cite[Chapter 8]{yaf}).
Then $\xi(\cdot;H,H_0)$ is called {\em the spectral shift
function} (SSF) for the operator pair $(H,H_0)$.\\
By the {\em Birman--Krein formula}, for almost every $E>0 = \inf \sigma_{\rm ac}(H)$, the SSF
$\xi(E; H, H_0)$ coincides with  {\em the scattering phase}
for the operator pair $(H,H_0)$ (see the original work \cite{bk} or the monograph \cite{yaf}).\\
Further, for almost every $E<0$ we have
$$
-\xi(E; H, H_0) = {\rm Tr}\,{\bf 1}_{(-\infty, E)}(H),
$$
 where ${\rm Tr}\,{\bf 1}_{(-\infty, E)}(H)$ is just the number of the eigenvalues of $H$ less than $E$, counted with their multiplicities.\\
The above properties follow directly from the general abstract theory of the SSF (see e.g. \cite[Chapter 8]{yaf}). By \cite[Proposition 2.5]{bpr}, the SSF for the operator pair defined in \eqref{1} and \eqref{2} possesses the following more particular features:
\begin{itemize}
\item  $\xi(\cdot; H,H_0)$ is bounded on every compact subset
of $\re \setminus 2b{\mathbb Z}_+$;

\item $\xi(\cdot; H,H_0)$ is continuous on $\re
\setminus (2b{\mathbb Z}_+ \cup \sigma_{\rm pp}(H))$ where
$\sigma_{\rm pp}(H)$ is the set of the eigenvalues of $H$.
\end{itemize}
Our first goal is to describe the asymptotic behavior of the SSF $\xi(E;H,H_0)$ as $E \to 2bq$, $q \in {\mathbb Z}_+$.
This behavior will be described in terms of auxiliary Berezin--Toeplitz operators studied in more detail in Subsection \ref{ss22}. The next subsection deals with the well-known properties of the Landau Hamiltonian, i.e. the 2D Schr\"odinger operator with constant magnetic fields.

\subsection{The Landau Hamiltonian}
\label{ss22a}
We have
    \bel{3}
H_0 = \hope  \otimes I_{\parallel} + I_{\perp}\otimes \hopa
    \ee
where $I_{\perp}$ and $I_{\parallel}$ are the identities in
$L^2(\rd_{\xp})$ and $L^2(\re_{x_3})$ respectively,
$$
 \hope : = \left(-i\frac{\partial}{\partial x_1} +
 \frac{bx_2}{2}\right)^2 + \left(-i\frac{\partial}{\partial x_2} -
 \frac{bx_1}{2}\right)^2 - b
$$
is the (shifted) {\em  Landau Hamiltonian}, self-adjoint in
$L^2(\rd_{\xp})$, and
$$
\hopa : = -\frac{d^2}{dx_3^2}
$$
is the 1D free Hamiltonian, self-adjoint in $L^2(\re_{x_3})$.
Note
that
$\hope  = a^* a$
where
$$
a: = -2i e^{-b|z|^2/4} \frac{\partial}{\partial \bar{z}} e^{b|z|^2/4}, \; \bar{z} = x_1
- i x_2,
$$
is {\em the magnetic annihilation operator}, and
 $$
a^*: = -2i e^{b|z|^2/4} \frac{\partial}{\partial z} e^{-b|z|^2/4}, \; z = x_1
+ i x_2,
$$
is {\em the magnetic creation operator}, adjoint to $a$ in $L^2(\rd)$.
Moreover,
$[a, a^*] = 2b$. Therefore, $\sigma(\hope) = \cup_{q=0}^{\infty}\{2bq\}$.
Furthermore,
$${\rm Ker}\,\hope = {\rm Ker}\,a = \Big\{f \in L^2(\rd) | f = g e^{-b|z|^2/4},
  \, \frac{\partial g}{\partial \bar{z}} = 0\Big\}
  $$
is the classical {\em Fock--Segal--Bargmann} space (see e.g.
\cite{ha}), and
$$
{\rm Ker}\,(\hope-2bq) = (a^*)^q {\rm Ker}\,\hope,
\quad q\geq
  1.
  $$
Evidently,
  $$
  {\rm dim}\,{\rm Ker}\,(\hope-2bq) = \infty
  $$
for each $q\in {\mathbb Z }_+$.\\
Representation \eqref{3} shows that the operator $H_0$ has {\em a waveguide structure} since the transversal operator $\hope$ has a purely point spectrum and the set of its eigenvalues is a discrete subset of the real axis, while the longitudinal operator $\hopa$ has a purely absolutely continuous spectrum. This waveguide structure of $H_0$ explains the qualification of the Landau levels, i.e. the eigenvalues of the transversal operator $\hope$, as {\em thresholds} in the spectrum of the ``total" 3D operator $H_0$. The important difference with the usual waveguides (see e.g. \cite{guil, chris}) is that the eigenvalues $2bq$, $q \in {\mathbb Z}_+$, of the transversal operator $\hope$ are of infinite multiplicity. We would like to underline here that most of the phenomena discussed in the present article are due to the infinite degeneracy of the Landau levels regarded as eigenvalues of $\hope$.

\subsection{Berezin--Toeplitz operators}
\label{ss22}

Fix $q \in {\mathbb Z}_+$. Denote by $p_q$ the orthogonal projection onto
${\rm Ker}\,(\hope-2bq)$. As discussed in the previous subsection, we have
${\rm rank}\; p_q = \infty$.\\
Let $U : \rd \to \ce$ be a Lebesgue measurable function.
Introduce {\em the Berezin--Toeplitz operator}
$$
p_q U p_q : p_q L^2(\rd)
\to p_q L^2(\rd).
$$
We will call $U$ {\em the symbol} of the operator $p_q U p_q$.
Evidently, if $U \in L^{\infty}(\rd)$ then $p_q U p_q$ is
bounded, and
$\|p_q U p_q\| \leq \|U\|_{L^{\infty}(\rd)}$.
Moreover, if $U \in L^{p}(\rd)$, $p \in [1,\infty)$, then by \cite[Lemma 5.1]{r0} or
\cite[Lemma 3.1]{fr},  we have $p_q U p_q \in S_p$,  the $p$th Schatten--von Neumann class, and
$$
\|p_q U p_q\|_{S_p}^p \leq \frac{b}{2\pi} \|U\|_{L^p(\rd)}^p.
$$
As a corollary, if $U \in L^1_{\rm loc}(\rd)$ and $U(x) \to 0$ as $|x| \to \infty$, then $p_q U p_q$ is compact. \\
Further, $p_0 U p_0$ with domain $p_0 L^2(\rd)$ is unitarily equivalent
to the $\Psi$DO : $L^2(\re) \to L^2(\re)$ with anti-Wick symbol
$$
\omega(y, \eta) : = U(b^{-1/2}\eta, b^{-1/2}y), \; (y,\eta) \in
T^*\re,
$$
while by \cite[Lemma 9.2]{bpr} the operator $p_q U p_q$ with any $q \in {\mathbb Z}_+$ is unitarily equivalent to
 \bel{5}
p_0 \left(\sum_{s=0}^q \frac{q!}{(2b)^s(s!)^2 (q-s)!} \Delta^s U\right)
p_0
    \ee
    which is quite useful when we want to reduce the analysis at the higher Landau levels to analysis at the first Landau level.
 Note that the differential operation occurring in \eqref{5} can be written as
${\rm L}_q\left(-\frac{\Delta}{2b}\right)$ where
$$
{\rm L}_q\left(t\right) : = \sum_{s=0}^q \frac{q!}{(s!)^2 (q-s)!} (-t)^s, \quad t\in \re,
$$
is the $q$th Laguerre polynomial. A more abstract point of view concerning the unitary equivalence between $p_q U p_q$ and the operator defined in \eqref{5} could be found in \cite{goff}.\\
The following three lemmas deal with the spectral asymptotics for compact Berezin--Toeplitz operators whose symbols $U$ admit respectively a power-like decay, an exponential decay, or have a compact support. More precisely, we discuss the asymptotics as $s \downarrow 0$ of the eigenvalue counting function ${\rm Tr}\,{\bf 1}_{(s,\infty)}(p_q U p_q)$.

\begin{lemma} \label{l1} {\rm \cite[Theorem 2.6]{r0}} Let $0 \leq U \in
C^1(\rd)$, and
$$
U(\xp) = u_0(\xp/|\xp|) |\xp|^{-\alpha}(1 + o(1)),
$$
$$
|\nabla U(\xp)| = O(|\xp|^{-\alpha-1}),
$$
as $|\xp| \to
\infty$, with $\alpha > 0$, and $0<u_0 \in C({\mathbb
S}^1)$.
Fix $q \in{\mathbb Z}_+$. Then
\bel{6}
{\rm Tr}\,{\bf 1}_{(s,\infty)}(p_q U p_q) =
\frac{b}{2\pi} \left|\left\{\xp \in \rd | U(\xp) > s\right\}\right| (1 + o(1)) =
\psi_{\alpha}(s) (1 + o(1)), \quad s \downarrow 0,
    \ee
where $|\cdot |$ denotes the Lebesgue measure, and
    \bel{san27}
\psi_{\alpha}(s) : = s^{-2/\alpha}
\frac{b}{4\pi}
\int_{{\mathbb S}^1} u_0(t)^{2/\alpha}dt.
  \ee
\end{lemma}

\begin{lemma} \label{l2} {\rm \cite[Theorem 2.1, Proposition 4.1]{rw}} Let
$0 \leq U \in L^{\infty}(\rd)$ and
$$
\ln{U(\xp)} = -\mu |\xp|^{2\beta} (1 + o(1)), \quad |\xp| \to
\infty,
$$
with $\beta \in (0,\infty)$,  $\mu \in (0,\infty)$.
Fix $q \in {\mathbb Z}_+$. Then
    \bel{7}
{\rm Tr}\,{\bf 1}_{(s,\infty)}(p_q U p_q) = \varphi_{\beta}(s)(1 + o(1)) , \quad s \downarrow
0,
    \ee
where
    \bel{san28}
\varphi_{\beta}(s) : = \left\{
\begin{array} {l}
\frac{b}{2\mu^{1/\beta}} |\ln{s}|^{1/\beta} \; {\rm if} \; 0 <
\beta < 1,\\
\frac{1}{\ln{(1+2\mu/b)}}|\ln{s}| \; {\rm if} \;
\beta = 1, \\
\frac{\beta}{\beta - 1}(\ln|\ln{s}|)^{-1}|\ln{s}| \; {\rm if} \;
1 < \beta < \infty.
\end{array}
\right.
    \ee
\end{lemma}

\begin{lemma} \label{l3} {{\rm \cite[Theorem 2.2, Proposition 4.1]{rw}}} Let
$0 \leq U \in L^{\infty}(\rd)$,  ${\rm supp}\,U$ be
compact, and  $U \geq C > 0$ on an open non-empty subset of
$\rd$. Fix $q \in {\mathbb Z}_+$.  Then
  \bel{8}
{\rm Tr}\,{\bf 1}_{(s,\infty)}(p_q U p_q) = \varphi_{\infty}(s) (1 + o(1)), \quad s \downarrow
0,
    \ee
where
    \bel{san29}
\varphi_{\infty}(s) : = (\ln|\ln{s}|)^{-1}|\ln{s}|.
    \ee
\end{lemma}
 Asymptotic relation \eqref{6} is of semiclassical nature in the sense that it is written in terms of the measure of that part of ``the phase space" $\rd$ where the symbol $U$ of the operator $p_q U p_q$ is greater than $s>0$. Similarly, asymptotic relation \eqref{7} with $\beta \in (0,1)$ is of semiclassical nature.
Asymptotic relation \eqref{7} with $\beta = 1$ is the border-line one: the order is semiclassical but the coefficient $\frac{1}{\ln{(1+2\mu/b)}}$ is not. Note that the main asymptotic term of $\frac{1}{\ln{(1+2\mu/b)}}$ as $b \to \infty$ coincides with the semiclassical coefficient $\frac{b}{2\mu}$. Finally, asymptotic relation \eqref{7} with $\beta \in (1,\infty)$ as well as asymptotic relation \eqref{8} are not of semiclassical nature.\\
Lemmas \ref{l1}, \ref{l2}, \ref{l3} have been cited in a similar form in several works of the present authors (see e.g. \cite{fr, bbr1, bbr2}). We have chosen this form since, in our opinion, it contains a reasonable scale of the possible types of decay of $V$, which, in particular, reveals clearly enough the passage from semiclassical to non semiclassical asymptotic behavior of the eigenvalue counting function for $p_qUp_q$. Of course, these three lemmas do not cover all possible symbols $U$ for which the main asymptotic term of ${\rm Tr}\,{\bf 1}_{(s,\infty)}(p_q U p_q)$ as $s \downarrow 0$ can be found explicitly. For example, if $U \geq C > 0$ on an open non-empty subset of
$\rd$, and
$$
\lim_{|\xp| \to \infty} \frac{\ln{(-\ln U(\xp))}}{\ln |\xp|} = \infty,
$$

then \eqref{7} and \eqref{8} easily imply that
$$
1 \leq \liminf_{s \downarrow 0} \frac{{\rm Tr}\,{\bf 1}_{(s,\infty)}(p_q U p_q)}{(\ln|\ln{s}|)^{-1}|\ln{s}|} \leq
\limsup_{s \downarrow 0} \frac{{\rm Tr}\,{\bf 1}_{(s,\infty)}(p_q U p_q)}{(\ln|\ln{s}|)^{-1}|\ln{s}|} \leq \frac{\beta}{\beta - 1}
$$
for any $\beta \gg 1$. Letting $\beta \to \infty$, we find that \eqref{8} again holds true.

\subsection{Asymptotics of $\xi(E; H,H_0)$ as $E \to 2bq$}
\label{ss24}

Let $V$ satisfy D. For $\xp \in \rd$, $\lambda \geq 0$, set
    \bel{san30}
W(\xp) : = \int_{\re} |V(\xp,x_3)| dx_3,
    \ee
$$
{\cal W}_{\lambda} = {\cal W}_{\lambda}(\xp) : = \left(
\begin{array} {cc}
w_{11} & w_{12}\\
w_{21} & w_{22}
\end{array}
\right),
$$
where
$$
w_{11} : = \int_{\re} |V(\xp, x_3)| \cos^2{(\sqrt{\lambda}x_3)} dx_3, \quad w_{22} : = \int_{\re} |V(\xp, x_3)|
\sin^2{(\sqrt{\lambda}x_3)} dx_3,
$$
$$
w_{12} = w_{21} : =
\int_{\re} |V(\xp, x_3)| \cos{(\sqrt{\lambda}x_3)}
\sin{(\sqrt{\lambda}x_3)}dx_3.
$$
Unless $V = 0$ almost everywhere, we have $${\rm rank}\,p_q W p_q = \infty, \quad {\rm rank}\,p_q
{\mathcal
  W}_{\lambda} p_q = \infty, \; \lambda \geq 0.
  $$
If $F_j(V; \lambda)$, $j=1,2$, are two real non decreasing functionals of $V$, depending on  $\lambda > 0$, we write
$$
F_1(V; \lambda) \sim F_2(V; \lambda), \quad \lambda \downarrow 0,
$$
if for each $\varepsilon \in (0,1)$ we have
$$
F_2((1-\varepsilon)V; \lambda) + O_\varepsilon(1) \leq
F_1(V; \lambda) \leq
 F_2((1+\varepsilon)V; \lambda) + O_\varepsilon(1).
$$
We also use analogous notations for non increasing functionals $F_j(V; \lambda)$ of $V$.

\begin{theorem} \label{th1} {\rm \cite[Theorems 3.1, 3.2]{fr}}
Let $V$ satisfy ${\rm D}_0$, and $V \geq 0$ or $V \leq 0$.
Fix
  $q \in {\mathbb Z}_+$. Then we have
    \bel{san1}
\xi(2bq-\lambda; H,H_0) = O(1), \quad \lambda \downarrow 0,
    \ee
if $V \geq 0$, and
    \bel{san2}
\xi(2bq-\lambda; H,H_0) \sim -
{\rm Tr}\,{\bf 1}_{(2\sqrt{\lambda},\infty)}\left(p_qWp_q\right),
\; \lambda \downarrow 0,
    \ee
if $V \leq 0$. Moreover,
    \bel{san3}
\xi(2bq+\lambda; H,H_0) \sim  \frac{1}{\pi}
{\rm Tr}\,\arctan{\left(\frac{p_q {\mathcal W}_{\lambda} p_q}{2\sqrt{\lambda}}\right)},
\; \lambda \downarrow 0,
    \ee
if $V \geq 0$, and
    \bel{san4}
\xi(2bq+\lambda; H,H_0) \sim  - \frac{1}{\pi}
{\rm Tr}\,\arctan{\left(\frac{p_q {\mathcal W}_{\lambda} p_q}{2\sqrt{\lambda}}\right)},
\; \lambda \downarrow 0,
    \ee
if $V \leq 0$.
\end{theorem}
Note that in the case $q = 0$ asymptotic relation \eqref{san2} concerns the distribution of the discrete eigenvalues of the operator $H$ with $V \leq 0$ near the origin which coincides with the infimum of its essential spectrum. Results of this type have been known for a long time, and could be found in:

\begin{itemize}
\item \cite{sob1, sob2, tam, r0, ivrii} in the case of {\em a power-like decay of $V$};

\item \cite{rw} in the case of {\em an exponential decay of $V$};

\item \cite{rw, melroz} in the case of {\em compactly supported potentials $V$}.
\end{itemize}
Inserting the results of Lemmas \ref{l1}, \ref{l2}, or \ref{l3} into \eqref{san2}, \eqref{san3}, and \eqref{san4}, we could obtain the main asymptotic term of the SSF as $E \to 2bq$.
We omit here these explicit formulae referring the reader to the original work (see \cite[Corollary 3.1]{fr}), and prefer to state here only the following intriguing

\begin{follow} \label{f1} {\rm \cite{r2}}
Let $V$ satisfy ${\rm D}_0$, and $V \leq 0$. Fix $q \in {\mathbb Z}_+$. Then
    \bel{san5}
\lim_{\lambda \downarrow 0} \frac{\xi(2bq + \lambda;
H,H_0)}{\xi(2bq - \lambda; H,H_0)}  = \frac{1}{2
\cos{\frac{\pi}{\alpha}}}
    \ee
if $W$ satisfies the assumptions of Lemma \ref{l1}, i.e. if $W$ admits a power-like decay with decay rate $\alpha > 2$, or
    \bel{san6}
\lim_{\lambda \downarrow 0} \frac{\xi(2bq + \lambda;
H,H_0)}{\xi(2bq - \lambda; H,H_0)}  = \frac{1}{2}
    \ee
if $W$ satisfies the assumptions of Lemma \ref{l2} or Lemma \ref{l3}, i.e. if $W$ decays exponentially\footnote{In the case of exponential decay of $W$ we should also suppose that $V$ satisfies ${\rm D}$ with $m_{\perp} > 2$ and $m_3>2$.} or has a compact support.
\end{follow}
 Relations \eqref{san5}--\eqref{san6} could be interpreted as {\em generalized Levinson formulae}. We recall that the classical Levinson formula relates
the number of the negative eigenvalues of $-\Delta + V$ with $V$ which decays sufficiently fast at infinity, and $\lim_{E \downarrow 0}\xi(E;-\Delta + V,-\Delta)$ (see the original work \cite{lev} or the survey article \cite{rob}).

\subsection{Sketch of the proof of Theorem \ref{th1}}
\label{ss26}
We start with a representation of the SSF due to A. Pushnitski \cite{p1, bpr}.
Assume that $V$ satisfies ${\rm D}$. Then the norm limit $$T(E) : = \lim_{\delta \downarrow 0} |V|^{1/2} (H_0 - E - i\delta)^{-1} |V|^{1/2}$$ exists for every $E \in \re \setminus 2b{\mathbb Z}_+$.
 Moreover, $T(E)$ is compact, and $0 \leq {\rm Im}\,T(E) \in S_1$ (see \cite[Lemma 4.2]{bpr}).
Assume in addition that $\pm V \geq 0$. Then for $E \in \re \setminus 2b{\mathbb Z}_+$ we have
 $$
 \xi(E; H, H_0)  =
 \pm \frac{1}{\pi} \int_\re {\rm Tr} {\bf 1}_{(1,\infty)}(\mp({\rm Re}\,T(E) + t {\rm Im}\,T(E))) \frac{dt}{1+t^2}
 $$
 (see \cite[Theorem 1.2]{p1}, \cite[Subsection 3.3]{bpr}).\\
The first important step in the proof of Theorem \ref{th1} is the estimate
    \bel{san7}
 \pm \xi(E; H, H_0) \sim
 \frac{1}{\pi} \int_\re {\rm Tr} {\bf 1}_{(1,\infty)}(\mp({\rm Re}\,T_q(E) + t {\rm Im}\,T_q(E))) \frac{dt}{1+t^2}, \quad
 E \to 2bq,
    \ee
  where
\begin{align*}
 T_q(E) : = {} & \lim_{\delta \downarrow 0} |V|^{1/2} (p_q \otimes I_{\parallel})(H_0 -E - i\delta)^{-1} |V|^{1/2} \\
= {} & \lim_{\delta \downarrow 0} |V|^{1/2} (p_q \otimes (\hopa + 2bq - E - i\delta)^{-1}) |V|^{1/2}, \quad E \neq 2bq.
\end{align*}
 If $E = 2bq - \lambda$ with $\lambda > 0$,
 then
 $T_q(E) = T_q(E)^*$,
 and \eqref{san7} implies
    \bel{san8}
 \pm \xi(E; H, H_0) \sim
 {\rm Tr} \, {\bf 1}_{(1,\infty)}(\mp\,T_q(E)), \quad
 E \to 2bq.
    \ee
 Moreover, we have  $T_q(E) \geq 0$, i.e. ${\rm Tr} \, {\bf 1}_{(1,\infty)}(-T_q(E)) = 0$. Then \eqref{san8} with the upper sign implies
 $$
 \xi(E; H, H_0) = O(1), \quad E \uparrow 2bq,
 $$
 provided that $V \geq 0$, i.e. we obtain \eqref{san1}.\\
Assume now that $V \leq 0$. The second important step in the proof of Theorem \ref{th1} is the estimate
    \bel{san10}
 {\rm Tr}\, {\bf 1}_{(1,\infty)}(T_q(2bq - \lambda)) \sim
 {\rm Tr}\, {\bf 1}_{(1,\infty)}\left(|V|^{1/2}\left(p_q \otimes S_-(\lambda)\right)|V|^{1/2}\right), \quad \lambda \downarrow 0,
    \ee
 where $S_-(\lambda)$ denotes the operator with constant integral kernel $\frac{1}{2\sqrt{\lambda}}$.
 Note that  $\frac{1}{2\sqrt{\lambda}}$ could be interpreted as the divergent part as $\lambda \downarrow 0$ of the integral kernel
 \bel{san17}
 \frac{e^{-\sqrt{\lambda}|x_3-x_3'|}}{2\sqrt{\lambda}}, \quad x_3, x_3' \in \re,
 \ee
 of the resolvent $(\hopa + \lambda)^{-1}$.
Our next step requires the following abstract
\begin{lemma} \label{l4} {\rm \cite[Theorem 8.1.4]{birsol}}
 Let $L$ be a linear compact operator acting between two, possible different, Hilbert spaces. Then for each $s>0$ we have
 $$
 {\rm Tr}\,{\bf 1}_{(s,\infty)}(L^* L) = {\rm Tr}\,{\bf 1}_{(s,\infty)}(L L^*).
 $$
 \end{lemma}
 Applying Lemma \ref{l4} with appropriate $L$, we immediately find that
    \bel{san11a}
    {\rm Tr}\, {\bf 1}_{(1,\infty)}\left(|V|^{1/2}\left(p_q \otimes S_-(\lambda)\right)|V|^{1/2}\right) =
    {\rm Tr}\, {\bf 1}_{(1,\infty)}\left(\frac{p_q W p_q}{2\sqrt{\lambda}}\right) = {\rm Tr}\, {\bf 1}_{(2\sqrt{\lambda},\infty)}\left(p_q W p_q\right).
    \ee
 Putting together \eqref{san8}, \eqref{san10},  and \eqref{san11a}, we obtain \eqref{san2}.\\
Let now $E = 2bq + \lambda$ with $\lambda \downarrow 0$. Then the next important step is the estimate
\begin{gather}
\frac{1}{\pi} \int_\re {\rm Tr} \, {\bf 1}_{(1,\infty)}(\mp({\rm Re}\,T_q(E) + t {\rm Im}\,T_q(E))) \frac{dt}{1+t^2} \sim
\frac{1}{\pi} \int_\re {\rm Tr} \, {\bf 1}_{(1,\infty)}(\mp t {\rm Im}\,T_q(E)) \frac{dt}{1+t^2} \nonumber  \\
= \frac{1}{\pi} {\rm Tr} \arctan{\left({\rm Im}\,T_q(E)\right)} = \frac{1}{\pi} {\rm Tr} \arctan{\left( |V|^{1/2} \left(p_q \otimes S_+(\lambda)\right) |V|^{1/2} \right)}    \label{san12a}
\end{gather}
where $S_+(\lambda)$ is the operator with integral kernel $\frac{\cos{\sqrt{\lambda}(x_3-x_3')}}{2\sqrt{\lambda}}$, $x_3, x_3' \in \re$. Applying Lemma \ref{l4} with appropriate $L$, we get
    \bel{san13a}
    \frac{1}{\pi} {\rm Tr} \arctan{\left( |V|^{1/2} \left(p_q \otimes S_+(\lambda)\right) |V|^{1/2} \right)} =
\frac{1}{\pi} {\rm Tr} \arctan{\left(\frac{p_q {\mathcal W}_{\lambda} p_q}{2\sqrt{\lambda}} \right)}.
    \ee
    Now the combination of \eqref{san7}, \eqref{san12a}, and \eqref{san13a}, yields \eqref{san3}--\eqref{san4}.

\subsection{Extensions of Theorem \ref{th1} to Pauli and Dirac operators}
\label{ss27}
Theorem \ref{th1} admits extensions to Pauli and Dirac operators with  non constant magnetic fields $(0,0,b)$ of constant direction. Here
$$
b = b_0 + \tilde{b},
$$
$b_0 \neq 0$ is a constant, and the function $\tilde{b} : \rd \to \re$ is such that the Poisson equation
$$
\Delta \tilde{\varphi} = \tilde{b}
$$
has a solution $ \tilde{\varphi} \in C_{\rm b}^2(\rd)$. In particular, $b$ may belong to a fairly large class of periodic or almost periodic functions of non zero mean value.\\
 In the case of the Pauli operator, the role of the Landau levels is played by the origin. The analogue of Theorem \ref{th1} could be found in \cite{r3}.
 Related results for negative energies (when the SSF is proportional to the eigenvalue counting function) are contained in \cite{it}.\\
In the case of the Dirac operator, the role of the Landau levels is played by the points $\pm m$ where $m >0$ is the mass of the relativistic quantum particle.
The analogue of Theorem \ref{th1} could be found in \cite{rtda}.

\section{Resonances near the Landau levels}
\label{s3}
\subsection{Embedded eigenvalues of $H$}
\label{ss31}
The singularity of the SSF as $E \uparrow 0$ in the case $V \leq 0$ has a simple explanation: the existence of infinitely many negative discrete eigenvalues of $H$ accumulating at the origin which coincides with the infimum of the essential spectrum of $H$. The explanation of the singularities of the SSF at the higher Landau levels is much less transparent. There is no evidence that in the general case these singularities are due (only) to the accumulation at the Landau levels of embedded eigenvalues of $H$. That is why the natural conjecture is that the singularities of the SSF at the higher Landau levels are related to the accumulation of resonances of $H$ at these levels; we discuss this possible accumulation in the several following subsections. \\
Our next theorem however stresses  the fact that the magnetic Hamiltonians in the presence of an appropriate symmetry are much apter to have embedded eigenvalues than the non magnetic ones.

\begin{theorem} \label{thsan1}
Let the operator $V(H_0 + 1)^{-1}$ be compact in  $L^2(\rt)$. Assume moreover, that
$V$ is axisymmetric, i.e. it depends only on $\rho : = |\xp|$ and $x_3$. \\
(i) Suppose that $V$ satisfies
 \bel{san11}
-2b < V(x) \leq -C {\bf 1}_K(x), \quad x \in \re^3,
    \ee
where $C>0$, and $K \subset \re^3$ is an open non empty set.
    Then each interval $$(2b(q-1), 2bq), \quad q\in {\mathbb N},$$ contains at least one (embedded) eigenvalue of $H$.\\
(ii) Suppose now that $V$  satisfies
     \bel{san12}
-2b < V(x) \leq -C {\bf 1}_{\tilde{K}}(\xp) \langle x_3\rangle^{-m_3}, \quad
    x = (\xp, x_3)\in \re^3,
    \ee
    where $C>0$, $m_3 \in (0,2)$, and ${\tilde{K}} \subset \re^2$ is an open non empty set.  Then each interval $$(2b(q-1), 2bq), \quad q\in {\mathbb N},$$ contains a sequence of (embedded) eigenvalues of $H$ which converges to $2bq$.
    \end{theorem}
    The first part of the theorem is contained in \cite[Theorem 5.1]{ahs}, and the simple modifications needed for the second part are briefly outlined in \cite[Subsection 3.1]{r2}. Nonetheless, due to some imprecise statements of the results of Theorem \ref{thsan1} which appeared in \cite[p. 385]{fr} and
    \cite[p. 3457]{br}, and were already commented in \cite[Subsection 3.1]{r2}, we include in the Appendix a detailed sketch of the proof of Theorem \ref{thsan1}.

On the contrary, it is expected that $H$ has no (embedded) eigenvalues in the case $V \geq 0$. The fact that the SSF rests bounded below each Landau level (see \eqref{san1}) confirms this conjecture. It has been proved for small $V$. More precisely, we have

\begin{theorem} {\rm \cite[Proposition 7]{bbr1}}
Let $V \geq 0$ satisfy ${\rm D}$ with $m_{\perp} > 0$ and $m_{3} > 2$. There exists $\kappa_{0} > 0$ such that, for any $0 \leq \kappa \leq \kappa_{0}$, $H_{0} + \kappa V$ has no (embedded) eigenvalues  in $\R \setminus 2 b \Z_{+}$.
\end{theorem}

Moreover, without the smallness assumption, \cite[Corollary 6.7]{bbr2} states that the (embedded) eigenvalues of $H$ form a discrete set for generic potentials $V \geq 0$ satisfying ${\rm D}_{\rm exp}$. Nevertheless, the absence of eigenvalues for general non-negative $V$ remains an open problem.

\subsection{Meromorphic continuation of the resolvent of $H$ and
definition of
  resonances}
  \label{ss32}
  As mentioned in the previous subsection, it is expected that in the generic case the singularities of the SSF described in Theorem \ref{th1} are related to the accumulation of the resonances
  of $H$ at the Landau levels. The first step in this investigation
  is, of course, the definition of the resonances themselves. As is
  generally accepted nowadays, we will define these resonances as the
  poles of a meromorphic extension of the  resolvent $(H-z)^{-1}$ to
  an appropriate Riemann surface ${\mathcal M}$.
For $z \in \ce_+ : = \{\zeta \in \ce \, | \, {\rm Im}\;\zeta >
0\}$ we have
$$
(H_0 - z)^{-1} = \sum_{q=0}^{\infty} p_q \otimes \left( \hopa + 2bq - z\right)^{-1}.
$$
Recall that the resolvent  $\left( \hopa
 - z\right)^{-1}$ with $z \in \ce_+$ admits the integral kernel
$$
-\frac{e^{i\sqrt{z} |x_3 - x_3'|}}{2i\sqrt{z}}, \quad x_3, x_3' \in \re, \quad {\rm Im}\,\sqrt{z} > 0,
$$
(cf. \eqref{san17}). Hence, for any $q \in {\mathbb Z}_+$ the operator $p_q \otimes ( \hopa + 2bq - z )^{-1}$ admits a standard analytic extension to the two-sheeted Riemann surface of the square root $\sqrt{z-2bq}$ which however depends on $q$. Therefore, we define ${\mathcal M}$
as  the infinite-sheeted Riemann surface of the countable family
    \bel{san18}
\left\{\sqrt{z-2bq}\right\}_{q \in {\mathbb Z}_+}.
    \ee
Let ${\mathcal P}_G : {\mathcal M} \to \ce\setminus 2b{\mathbb Z}_+$
be the corresponding covering.\\
The properties of the Riemann surface ${\mathcal M}$ have been studied in detail in \cite[Section 2]{bbr1}.
Similar infinite-sheeted Riemann surfaces appearing in the spectral and resonance theory for perturbed waveguides have been introduced e.g. in
\cite{guil, edw, chris}. In this case the family analogous to \eqref{san18} is
$$
\left\{\sqrt{z-\mu_q}\right\}_{q \in \Z_{+}}
$$
where $\mu_q$, $q \in \Z_{+}$, are the distinct eigenvalues of the transversal operator which in the case of a  waveguide is lower bounded, and has a discrete spectrum.\\
The global structure of the Riemann surface ${\mathcal M}$ is quite complicated and may make difficult the analysis of the resonances of $H$. The investigation of their asymptotic distribution near a fixed Landau level $2bq$, $q \in {\mathbb Z}_+$, however is facilitated by the fact that in this case we are concerned with the local properties of  ${\mathcal M}$, and in a domain analytically diffeomorphic to a vicinity of $2bq$, the surface ${\mathcal M}$ resembles the two-sheeted Riemann surface  of the square root $\sqrt{z-2bq}$. Namely, if we put
$$
D(\lambda_0, \varepsilon) : = \{\lambda \in {\mathbb C} | \,|\lambda
  - \lambda_0| < \varepsilon\}, \quad
D(\lambda_0, \varepsilon)^{*}: =
  \{\lambda \in {\mathbb C} | \,0<|\lambda
  - \lambda_0| < \varepsilon\},
$$
for $\lambda_0 \in {\mathbb C}$ and $\varepsilon
  > 0$, then there exists a domain $D_q^* \subset {\mathcal M}$, and a
an analytic bijection
    \bel{san20}
 D (0 , \sqrt{2b})^{*} \ni k \mapsto z_{q} (k) \in
D_{q}^{*} \subset {\mathcal M},
 \ee
such that ${\mathcal P}_{G} (z_{q} (k)) = 2 b q + k^{2}$.\\
For $N>0$ denote by
${\mathcal M}_{N}$ the part of ${\mathcal M}$ where
${\rm Im}\, \sqrt{z-2b q} > -N$ for all $q \in {\mathbb Z}_+$.
Then, $\cup_{N>0} {\mathcal M}_{N} = {\mathcal M}$.

\begin{pr} \label{pr1} {\rm \cite[Propositions 1,2]{bbr1}}
{\rm (i)} For each $N > 0$ the operator-valued function
$$
(H_0-z)^{-1} : e^{-N \langle x_3
\rangle} L^2 (\rt ) \rightarrow e^{N \langle x_3
\rangle} L^2 (\rt)
$$
has an analytic  extension from
$\ce_+$ to  ${\mathcal M}_{N}$.\\
{\rm (ii)} Suppose that $V$ satisfies ${\rm D}_{\rm exp}$ with $\mper > 0$. Then
  for each $N > 0$ the
operator-valued function
$$
(H-z)^{-1} : e^{-N \langle x_3 \rangle} L^2 (\rt) \rightarrow  e^{N \langle x_3 \rangle} L^2 (\rt),
$$
has a meromorphic extension from $\ce_+$ to ${\mathcal
  M}_{N}$ whose poles and  residue ranks  do not depend on $N$.
  \end{pr}
We define the resonances of $H$ as the poles of the meromorphic extension of
the resolvent $(H-z)^{-1}$, and denote their set by $\res (H)$. For $z_0 \in \res (H)$ define its multiplicity by
$$
{\rm mult}\,(z_0): = {\rm rank}\; \frac{1}{2i\pi} \int_{\gamma} (H-z)^{-1} dz,
$$
where $\gamma$ is a circle centered at
$z_0$ and run over in the clockwise direction, such that $\overline{{\rm Int}\,\gamma}$ contains no elements of $\res (H)\setminus\{z_0\}$.

\subsection{Resonance-free regions and regions with infinitely many resonances}
\label{ss33}
One of the main technical achievements of our article \cite{bbr1} was
the identification of the resonances (together with their
multiplicities) as the zeroes of an appropriate
2-determinant.
We recall  that for a Hilbert--Schmidt operator $T$ the 2-determinant is defined as
$$
{\rm det}_2(I+T) = {\rm det}(I+T)e^{-T}.
$$

\begin{pr} \label{pr2} {\rm \cite[Proposition 3]{bbr1}}
Suppose that $V$ satisfies ${\rm D}_{\rm exp}$ with $\mper > 2$.
Introduce ${\mathcal T}_V(z)$, the analytic extension from $\ce_+$ to ${\mathcal
  M}_{N}$, of
 \bel{san50}
{\mathcal T}_V(z) : = {\rm sign}\,{V}\, \vert V \vert^{1/2} (H_0 - z)^{-1} \vert V \vert^{1/2}.
    \ee
Then
$z_0 \in {\mathcal M}$ is a resonance of $H$ if and only if
$-1$ is an eigenvalue of ${\mathcal T}_V(z_{0})$. Moreover,
    \bel{san21}
{\det}_2 \big( (H-z)(H_0-z)^{-1} \big) = {\det}_2 \big( I + {\mathcal T}_V(z) \big)
    \ee
has an analytic continuation from ${\mathbb C}_{+}$ to ${\mathcal M}$ whose zeroes are the resonances of $H$, and if $z_0$ is a resonance, then there exists a holomorphic function $f (z)$, for $z$ close to $z_0$, such that $f(z_0) \neq 0$ and
$$
{\det}_2 \big( I + {\mathcal T}_V(z) \big) = (z-z_0)^{{\rm mult}\,{(z_0)}} f ( z ) .
$$
\end{pr}

In the case of a trace-class perturbation $H-H_0$ the determinant ${\det} \big( (H-z)(H_0-z)^{-1} \big)$, $z \in \C_+$, coincides with the classical perturbation determinant introduced by M. G. Krein in \cite{krein} (see also \cite[Section IV.3]{gk}). In the case of Hilbert-Schmidt perturbations $H - H_0 \in S_2$ (or relatively Hilbert--Schmidt perturbations $(H - H_0)(H_0 - z)^{-1} \in S_2$) the (generalized) perturbation 2-determinants were introduced by L. S. Koplienko in \cite{kop} where he considered as well the whole Schatten--von Neumann scale $H - H_0 \in S_r$ (or  $(H - H_0)(H_0 - z)^{-1} \in S_r$), $r \geq 1$. \\
\begin{figure}
\begin{center}
\begin{picture}(0,0)%
\includegraphics{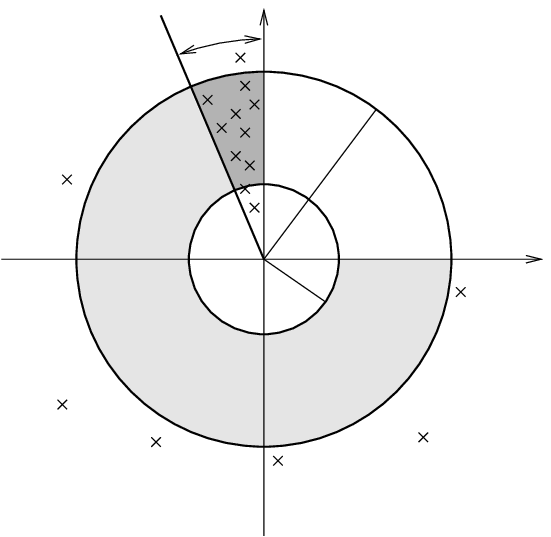}%
\end{picture}%
\setlength{\unitlength}{1184sp}%
\begingroup\makeatletter\ifx\SetFigFont\undefined%
\gdef\SetFigFont#1#2#3#4#5{%
  \reset@font\fontsize{#1}{#2pt}%
  \fontfamily{#3}\fontseries{#4}\fontshape{#5}%
  \selectfont}%
\fi\endgroup%
\begin{picture}(8744,8519)(1779,-9608)
\put(6901,-3211){\makebox(0,0)[lb]{\smash{{\SetFigFont{9}{10.8}{\rmdefault}{\mddefault}{\updefault}$r_{0}$}}}}
\put(10351,-4936){\makebox(0,0)[lb]{\smash{{\SetFigFont{9}{10.8}{\rmdefault}{\mddefault}{\updefault}$\operatorname{Re} k$}}}}
\put(6376,-1336){\makebox(0,0)[lb]{\smash{{\SetFigFont{9}{10.8}{\rmdefault}{\mddefault}{\updefault}$\im k$}}}}
\put(4951,-1411){\makebox(0,0)[lb]{\smash{{\SetFigFont{9}{10.8}{\rmdefault}{\mddefault}{\updefault}$\delta$}}}}
\put(6601,-5461){\makebox(0,0)[lb]{\smash{{\SetFigFont{9}{10.8}{\rmdefault}{\mddefault}{\updefault}$r$}}}}
\put(8926,-2086){\makebox(0,0)[lb]{\smash{{\SetFigFont{9}{10.8}{\rmdefault}{\mddefault}{\updefault}$V \leq 0$}}}}
\end{picture} $\qquad \qquad$ \begin{picture}(0,0)%
\includegraphics{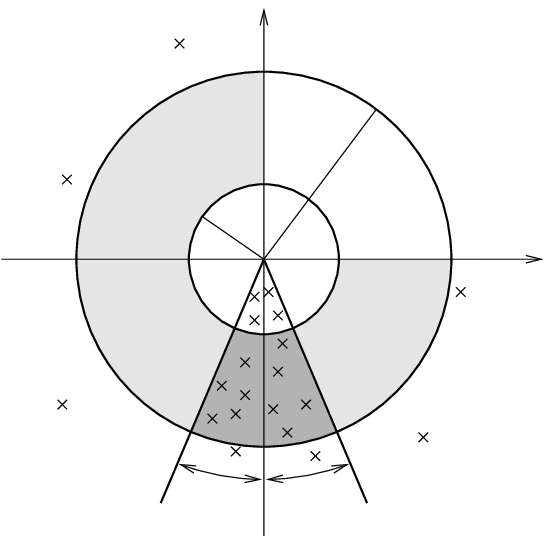}%
\end{picture}%
\setlength{\unitlength}{1184sp}%
\begingroup\makeatletter\ifx\SetFigFont\undefined%
\gdef\SetFigFont#1#2#3#4#5{%
  \reset@font\fontsize{#1}{#2pt}%
  \fontfamily{#3}\fontseries{#4}\fontshape{#5}%
  \selectfont}%
\fi\endgroup%
\begin{picture}(8744,8519)(1779,-9608)
\put(4951,-9061){\makebox(0,0)[lb]{\smash{{\SetFigFont{9}{10.8}{\rmdefault}{\mddefault}{\updefault}$\delta$}}}}
\put(6901,-3211){\makebox(0,0)[lb]{\smash{{\SetFigFont{9}{10.8}{\rmdefault}{\mddefault}{\updefault}$r_{0}$}}}}
\put(6526,-9061){\makebox(0,0)[lb]{\smash{{\SetFigFont{9}{10.8}{\rmdefault}{\mddefault}{\updefault}$\delta$}}}}
\put(5476,-4711){\makebox(0,0)[lb]{\smash{{\SetFigFont{9}{10.8}{\rmdefault}{\mddefault}{\updefault}$r$}}}}
\put(10351,-4936){\makebox(0,0)[lb]{\smash{{\SetFigFont{9}{10.8}{\rmdefault}{\mddefault}{\updefault}$\operatorname{Re} k$}}}}
\put(8926,-2086){\makebox(0,0)[lb]{\smash{{\SetFigFont{9}{10.8}{\rmdefault}{\mddefault}{\updefault}$V \geq 0$}}}}
\put(6376,-1336){\makebox(0,0)[lb]{\smash{{\SetFigFont{9}{10.8}{\rmdefault}{\mddefault}{\updefault}$\im k$}}}}
\end{picture}%
\end{center}
\caption{Resonances near a Landau level for $V$ of definite sign, concentrated near the semi-axis $k =- i (\hbox{sgn} V)  (0, + \infty)$.} \label{figres}
\end{figure}

The identification of the resonances of $H$ as zeroes of the
2-determinant \eqref{san21} allowed us to obtain in \cite{bbr1}
various local estimates of the number of resonances near any fixed
Landau level $2bq$, $q \in {\mathbb Z}_+$. As a typical and important
example, we state below Theorem \ref{th2}. It shows that if $V$ is
small, and its sign is fixed, then the
possible resonances of $H$ near $2bq$, parametrized according to
\eqref{san20}, are concentrated in sectors centered at $2bq$,
adjoining the imaginary axis, and depending on the sign of
$V$, as illustrated in Figure \ref{figres}. Moreover,
for rapidly decaying $V$, the number of the resonances of $H$ in each of these sectors is infinite.

\begin{theorem} \label{th2}  {\rm \cite[Theorem 2]{bbr1}} Let $0 < r_0 <\sqrt{2b}$ and $q \in {\mathbb Z}_+$.
Assume $V$ satisfies ${\rm
 D}_{\rm exp}$  with $\mper >2$, and is of definite sign $J$. Then
for any $\delta >0$ there exists $\varkappa_0 >0$ such that:\\
{\rm (i)} $H_0 + \varkappa
V$ has no resonances in $$\{z = z_q(k) \, | \, 0< |k |< r_0 , \; -J
{\rm Im}\; k \leq \frac{1}{\delta} |{\rm Re}\; k|\}$$
 for any $0 \leq \varkappa \leq \varkappa_0$.\\
{\rm (ii)}  If the function $W$ defined in \eqref{san30}, satisfies $\ln W(\xp) \leq - C\langle \xp
\rangle^2$, then  for any $0 < \varkappa \leq \varkappa_0$, the operator $H_0 +
\varkappa V$ has  an infinite number of resonances in
$$
\{z = z_q(k) \, | \,0< |k |< r_0 , \; -J {\rm Im}\; k > \frac{1}{\delta}
|{\rm Re}\; k|\}.
$$
\end{theorem}
In the above result, the assumption $\mper >2$ is not
  necessary. It was made in \cite{bbr1} in order to define the
  2-determinant, but using the notion of index (see Subsection
  \ref{ss35a}), Theorem \ref{th2} can be proved also for $\mper >0$.

\subsection{Asymptotics of the resonance counting function}
\label{ss35}
In spite of the undisputable usefulness of the methods developed and applied in \cite{bbr1}, they led to a loss of sharpness in some of the crucial estimates which hindered us to obtain the main asymptotic term of the number of the resonances of $H$ lying on an annulus centered at the Landau level $2bq$, $q \in {\mathbb Z}_+$, as the inner radius of the annulus tends to zero. This result central in the theory of the resonances of $H$, was obtained recently in \cite{bbr2} using methods different from those of \cite{bbr1}, and is contained below in Theorem \ref{th3}.\\
For $q \in {\mathbb Z}_+$ and $z \in D(0, \sqrt{2b})$ set
$$
{\mathcal A}_q(z) : = J|V|^{1/2}\Big(\frac{1}{2} \Big(p_q \otimes e^{z|x_3-x_3'|}\Big) -
z  \sum_{j \neq q}  \Big(p_j \otimes \Big(\hopa + 2b(j-q) + z^2\Big)^{-1}\Big)\Big) |V|^{1/2}.
$$
Note that for $q \in \N$ and $k \in D(0,\sqrt{2b})^*$ we have
    \bel{san22}
I + {\mathcal T}_V(z_q(k)) = I - \frac{{\mathcal A}_q(ik)}{ik},
    \ee
    the operator-valued function ${\mathcal T}_V(z)$ being defined in \eqref{san50}.
    Let $\Pi_q$ be the orthogonal projection onto ${\rm Ker}\, {\mathcal A}_q(0).$
\begin{theorem} \label{th3}  {\rm \cite[Theorem 6.5]{bbr2}} Let $V$
  satisfy ${\rm D}_{\rm exp}$ with $\mper > 0$ and have a definite sign $J = \pm 1$. Let
the function $W$ defined in \eqref{san30}
satisfy the assumptions of Lemma \ref{l1}, \ref{l2}, or \ref{l3}. Fix $q \in {\mathbb Z}_+$, and assume that  $I - {\mathcal A}_q'(0) \Pi_q$ is invertible.
  Then the conclusions of Theorem \ref{th2} hold for $r_0$ small
    enough and we have
$$
\sum_{z_q(k) \in \res (H):\\r < |k| < r_0} {\rm mult}\,(z_q(k)) =
{\rm Tr}\, {\bf 1}_{(2r, \infty)}(p_q W p_q) ( 1 + o(1) )
$$
as $r \downarrow
0$.
\end{theorem}
A sketch of the proof of Theorem \ref{th3} can be found in the next subsection. Here we make some brief comments on its hypotheses as well on some related results. \\
First, the assumption that $I - {\mathcal A}_q'(0) \Pi_q$ is
invertible is  essential
and does not have a purely technical character. As shown in
\cite[Section 7]{bbr2}, abstract results close by spirit to Theorem
\ref{th3} may cease to be valid if this assumption is canceled. On the
other hand, it holds true for generic potentials $V$ because it
  is satisfied for $gV$ provided that $g \in \R$ is not in a discrete
  set ($1/g$ has to be distinct from the eigenvalues of the compact operator $ {\mathcal A}_q'(0) \Pi_q$). A simple sufficient condition is that the norm $\|V\|_{L^{\infty}(\re^3)}$ is small enough. \\
As explained at the end of Subsection \ref{ss22} the assumptions of Lemma \ref{l1}, \ref{l2}, or \ref{l3} could be replaced by more general ones. We stick to these assumptions of the sake of a coherent exposition of the article. \\
Theorem \ref{th3} has been extended recently by Sambou in \cite{sambou} to the setting of Pauli and Dirac Hamiltonians with non constant magnetic fields. In the case of the Pauli (resp., Dirac) operator, Sambou obtained the main asymptotic term of the resonance counting function for an annulus centered at the origin (resp., for annuli centered at $\pm m$). The class of the magnetic fields considered in  \cite{sambou} is quite close to the one described in Subsection \ref{ss27} above.\\
In a more general context, Theorem \ref{th3} belongs to a large group of results concerning the asymptotic behavior of various resonance counting functions. Among the best known problems of this type, notorious for its hardness, is the problem of finding the first asymptotic term as $r \to \infty$ of the number $N(r)$ of the resonances lying on the disk $\{z \in \C \, | \, |z| < r\}$ for the operator $-\Delta + V$ with, say, compactly supported $V$, self-adjoint in $L^2(\re^n)$, $n \geq 1$. Actually, the first asymptotic term as $r \to \infty$ of $N(r)$ is known in the general case only if $n = 1$. Then we have
    \bel{san23}
    N(r) = \frac{2a}{\pi} r(1 + o(1)), \quad r \to \infty,
    \ee
    where $a$ is the diameter (not the Lebesgue measure!) of the
    support of $V$ (see \cite{zw1, froese}). Hence, asymptotic
    relation \eqref{san23} is not of Weyl type. In the case $n \geq
    3$, $n$ odd, only sharp upper bounds of $N(r)$ are known in the
    general case (see e.g. \cite{zw2, vodev, st}). The main asymptotic
    term of $N(r)$ is known only in the exceptional cases when $V$ is
    radial and satisfies some additional properties (see \cite{zw2, st}); again this main term is not of Weyl type. Note that lower bounds for generic perturbations
      are also known (see for instance \cite{chris2}) and recently
      Sj\"ostrand \cite{sjostrand} obtained a probabilistic Weyl law for
      random perturbations. \\
      For more
    information on the asymptotics of $N(r)$ as $r \to \infty$, we
    refer the reader to the evolving lecture notes
    \cite{zw3}.

\subsection{Sketch of the proof of Theorem \ref{th3}}\label{ss35a}
In order to outline the proof of Theorem \ref{th3}, we need the following
abstract results.
Let $\D$ be a domain of $\C$ containing $0$, and let $\cH$ be a separable Hilbert space. Consider the analytic  function
\begin{equation*}
A: \D \longrightarrow S_{\infty} ( \cH ).
\end{equation*}
Let $\Pi(A)$ be the orthogonal projection onto ${\rm Ker}\,A(0)$.\\
In the sequel we will suppose that the following assumptions are fulfilled:
\begin{itemize}
\item ${\mathcal C_1}$: The operator $A(0)$ is self-adjoint;
\item ${\mathcal C_2}$: The operator $I - A'(0) \Pi(A)$ is invertible.
\end{itemize}
Let $\Omega \subset \D \setminus \{0\}$. Define the characteristic values of $I - A(z)/z$ on $\Omega$ as the points $z \in \Omega$ for which the operator $I - A(z)/z$ is not invertible. We will denote the characteristic values of $I - A(z)/z$ on $\Omega$ by $\cZ_A ( \Omega )$.
By ${\mathcal C_1}$ and ${\mathcal C_2}$ the set $\cZ_A ( \Omega )$ is
discrete. For $z_0 \in \cZ_A ( \Omega )$ define its multiplicity
 by
$$
{\rm Mult}(z_0) : =
\frac{1}{2\pi i} {\rm Tr}\,\int_{\gamma} \Big(I-\frac{A (z)}{z}\Big)'\Big(I-\frac{A (z)}{z}\Big)^{-1} dz
$$
where $\gamma$ is an appropriate circle centered at $z_0$.
\begin{pr} \label{pr3} {\rm \cite[Corollary 3.4]{bbr2}}
Assume ${\mathcal C}_1$ and ${\mathcal C}_2$. Suppose that the origin is an accumulation point of $\cZ_A(\D\setminus\{0\})$.
Then we have
$$|{\rm Im}\,z_0| = o(|z_0|), \quad  z_0 \in \cZ_A(\D\setminus\{0\}),
$$
as $z_0 \to 0$.  If, moreover, $\pm A(0) \geq 0$, then
$\pm {\rm Re}\,z_0 \geq 0$
for $z_0 \in \cZ_A(\D\setminus\{0\})$ with $|z_0|$ small enough.
\end{pr}
 Set
\begin{equation*}
\cN_A ( \Omega ):= \sum_{z_0 \in \cZ_A ( \Omega )} {\rm Mult}(z_0).
\end{equation*}
If
$\partial \Omega$ is sufficiently regular, and $\cZ_A ( \Omega ) \cap \partial \Omega = \emptyset$,  then we have
$$
\cN_A ( \Omega ) = {\rm ind}_{\partial \Omega} \Big(I-\frac{A (z)}{z}\Big) : =
\frac{1}{2\pi i} {\rm Tr}\,\int_{\partial \Omega} \Big(I-\frac{A (z)}{z}\Big)'\Big(I-\frac{A (z)}{z}\Big)^{-1} dz.
$$
The index ${\rm ind}_{\partial \Omega} \Big(I-\frac{A (z)}{z}\Big)$ plays a central role in the proof of Theorem \ref{th3}. More information about its properties could be found in
\cite{gs}, \cite[Section 4]{gl}, and \cite[Section 2]{bbr2} (see also
\cite{sjostrand} where the notion of index allows to define generalized determinants). \\
For $0 < a < b < \infty$ and $\theta > 0$ define the domain
\begin{equation} \label{a1}
C_\theta(a,b) : = \left\{x+iy \in \C \, | \, a<x<b, \quad |y| < \theta x\right\}.
\end{equation}
\begin{pr} \label{pr4} {\rm \cite[Corollary 3.11]{bbr2}}  Assume ${\mathcal C}_1$ and ${\mathcal C}_2$. Suppose moreover that
$$
{\rm Tr}\,{\bf 1}_{(r, \infty)}(A(0)) = \Phi(r) (1 + o(1)), \quad r \downarrow 0,
$$
where $\Phi$ satisfies $\Phi (r) \to \infty$ as $r \downarrow 0$, and
    \bel{san26}
    \Phi(r(1\pm\delta)) = \Phi(r)(1 + o(1) + O(\delta)), \quad r \downarrow 0,
    \ee
    for each sufficiently small $\delta>0$.
  Then  we have
$$
{\mathcal N}_A(C_\theta(r,1)) =  \Phi(r) (1 + o(1)), \quad r \downarrow 0,
$$
for any $\theta > 0$.
\end{pr}
It is easy to check that the functions $\Phi(r) = C r^{-\gamma}$, $\Phi(r) = C |\ln{r}|^{\gamma}$, or $\Phi(r) = C \frac{|\ln{r}|}{\ln{|\ln{r}|}}$, with some $\gamma , C > 0$, satisfy asymptotic   relation \eqref{san26}. Hence the functions $\psi_\alpha$, $\varphi_\beta$, and $\varphi_\infty$ defined respectively in \eqref{san27}, \eqref{san28}, and \eqref{san29}, satisfy it as well. \\
Now we are in position to prove Theorem \ref{th3}.
By \eqref{san22} we have that $z_q(k) \in \res (H)$ if and only if
$ik$ is a characteristic value of $I- {\mathcal A}_q(z) / z$.
Moreover,
$$
{\rm mult}\,(z_q(k)) = {\rm Mult}\,(ik).
$$
By Proposition \ref{pr3}  with $A = {\mathcal A}_q$,
$$
\left\{z_q(k) \in \res (H) \, | \, r < |k| < r_0\right\} =
\left\{z_q(k) \in \res (H) \, | \,\pm ik \in C_\theta(r,r_0) \right\}
+\, O(1), \quad r \downarrow 0.
$$
Now the claim of Theorem \ref{th3} follows from Proposition \ref{pr4}
with $A = {\mathcal A}_q$ combined with  Lemmas \ref{l1}, \ref{l2}, and \ref{l3}, since, by Lemma \ref{l4} with appropriate $L$, we have
${\rm Tr}\,{\bf 1}_{(r, \infty)}({\mathcal A}_q(0)) = {\rm Tr}\,{\bf 1}_{(2r, \infty)}(p_q W p_q)$.

\Subsection{Link between the SSF and the resonances}
The spectral shift function and the resonances are usually connected by the Breit--Wigner formula. This formula represents the derivative of the SSF as a sum of a harmonic measure associated to the resonances, and the imaginary part of a holomorphic function. When the resonances are close to the real axis, the Breit--Wigner approximation can be exploited to obtain asymptotic expansions of the SSF (see e.g. \cite{BrPe03_01}). On the contrary, it can be used to localize  resonances as in \cite[Appendix]{NaStZw03_01}. Eventually, it can imply local trace formulas in the spirit of \cite{Sj97_01}.\\

\begin{theorem}\label{a2} {\rm \cite[Theorem 3]{bbr1}} Let $V$ satisfies ${\rm D}_{\rm exp}$ with $\mper >2$. Then for $q \in \Z_{+}$ and $\varepsilon , \theta > 0$, there exist $r_{0} > 0$ and functions $g_{\pm} ( \cdot , r)$ holomorphic in  $\pm C_{\theta} ( 1 , 2 )$, such that, for $\lambda \in r [ 1 + \varepsilon , 2 -\varepsilon ]$, we have
\begin{equation*}
\xi ' ( 2 b q \pm \lambda ; H , H_{0} ) = \sum_{\fract{2 b q \pm w \in \res (H)}{w \in r C_{\theta} ( 1 , 2 ) \setminus \R}} \frac{\im w}{\pi |\lambda - w|^2} - \sum_{\fract{2 b q \pm w \in \res (H)}{w \in r [1 , 2 ]}} \delta ( \lambda - w ) +  \frac{1}{r} \im g_{\pm}' \Big( \frac{\lambda}{r} , r \Big) ,
\end{equation*}
where  $g_\pm(z, r)$ satisfies the estimate
\begin{equation*}
g_{\pm} ( z , r ) = O \big( | \ln r| \; r^{- \frac{1}{\mper}} \big) ,
\end{equation*}
uniformly with respect to $0 < r < r_{0} $ and $z \in C_{\theta} ( 1 + \varepsilon , 2 - \varepsilon )$.
\end{theorem}

A more general statement and some applications of this formula can be found in \cite{bbr1}. Note that this Breit--Wigner approximation implies that the SSF is analytic outside of the resonances (including the embedded eigenvalues) and their complex conjugate.\\
We close this subsection with the remark that we are not able yet to make the link between Theorem \ref{th1} and Theorem \ref{th3} using Theorem \ref{a2} for $V$ of definite sign. We believe that such a result would shed additional light on the relation between the SSF threshold singularities and the resonances of $H$, and hope to obtain it in a future work. Nevertheless, we can formally deduce the localization of the resonances from the asymptotics of the SSF. Indeed, for $V \leq 0$, the singularities of the main term in  \eqref{san2} are the numbers $2 b q + k^{2}$ with $2 k \in i \sigma ( p_{q} W p_{q} )$, which is in agreement with Theorem \ref{th3}. For $V \geq 0$, the uniform bound in \eqref{san1} should be explained by some cancelations  in the Breit--Wigner formula due to the particular form of the Riemann surface ${\mathcal M}$ and the localization of the resonances. Finally, the leading terms of \eqref{san3} and \eqref{san4} have singularities at $2 b q + k^{2}$ with  $2 k \in \pm i \sigma ( p_{q} W p_{q} )$ which can be explained by Theorem \ref{th3} and the symmetry of the singularities of the SSF with respect to the real axis.

\section{Appendix: Sketch of the proof of Theorem \ref{thsan1}}
\label{app}
Passing to cylindrical coordinates $(\rho, \varphi, x_3)$, and decomposing $u \in {\rm Dom}\,H_0$ into a Fourier series with respect to $\varphi$, we find that the operator $H_0$ is unitarily equivalent to the orthogonal sum $\sum_{m \in {\mathbb Z}} \oplus H_0^{(m)}$ where
$$
H_0^{(m)} : = -\frac{1}{\rho} \frac{\partial }{\partial \rho} \rho \frac{\partial }{\partial \rho} + \left(b\rho + \frac{m}{\rho}\right)^2 - \frac{\partial^2 }{\partial x_3^2} - b
$$
is self-adjoint in $L^2(\re_+ \times \re; \rho d\rho dx_3)$, while $H$  is unitarily equivalent to  $\sum_{m \in {\mathbb Z}} \oplus H^{(m)}$ with $H^{(m)} : = H_0^{(m)} + V(\rho, x_3)$. We have
$$
\sigma(H_0^{(m)}) = [2bm_+, \infty), \quad m \in {\mathbb Z}_+,
$$
where $m_+ : = \max\{m,0\}$ is the positive part of $m$.  Fix $q \in {\mathbb N}$. Since the operator $V(H_0 + 1)^{-1}$ is compact in $L^2(\rt)$, the operator  $V(H_0^{(q)} + 1)^{-1}$ is compact in $L^2(\re_+ \times \re; \rho d\rho dx_3)$. Therefore, the eigenvalues of the operator $H^{(q)}$ lying on the interval $(2b(q-1), 2bq)$ are discrete, and at the same time they are embedded eigenvalues of the ``total" operator $H$ since $\sigma_{{\rm ess}}(H) = [0,\infty)$. Let us estimate from below the quantity ${\rm Tr}\,{\bf 1}_{(2b(q-1), 2bq)}(H^{(q)})$. By the lower bounds in \eqref{san11}--\eqref{san12} we have
$$
{\rm Tr}\,{\bf 1}_{(2b(q-1), 2bq)}(H^{(q)}) = {\rm Tr}\,{\bf 1}_{(-\infty, 2bq)}(H^{(q)}).
$$
Evidently, it suffices to show that
    \bel{san15}
    {\rm Tr}\,{\bf 1}_{(-\infty, 2bq)}(H^{(q)}) \geq 1
    \ee
    in order to prove the first part of Theorem \ref{thsan1}. Moreover, since the spectrum of $H^{(q)}$ on $(-\infty, 2bq)$ is discrete and can accumulate only at $2bq$, it suffices to show that
    \bel{san16}
    {\rm Tr}\,{\bf 1}_{(-\infty, 2bq)}(H^{(q)}) = \infty
    \ee
    in order to prove the second part of this theorem. Let $\phi_q : \re_+ \to \re$ be an eigenfunction of the transversal part of the operator  $H^{(q)}$ satisfying
    $$
    -\frac{1}{\rho} \frac{d }{d \rho} \left(\rho \frac{d \phi_q }{d \rho}\right) + \left(b\rho + \frac{q}{\rho}\right)^2 \phi_q - b\phi_q = 2bq \phi_q,
    $$
    and $\int_0^{\infty} \phi_q^2 \rho d\rho = 1$. On  ${\rm H}^2(\re_{x_3})$ introduce the 1D Schr\"odinger operator $h_q : = - \frac{d^2 }{d x_3^2} + v_q$
    where
    $$
    v_q(x_3) : = \int_0^{\infty} V(\rho,x_3) \phi_q(\rho)^2 \rho d\rho, \quad x_3 \in \re.
    $$
    Restricting the quadratic form of the operator $H^{(q)}$ onto the subspace
    $$
    \left\{u : \re_+ \times \re \to {\mathbb C}\,|\,u(\rho,x_3) = \phi_q(\rho)\,w(x_3), \; w \in {\rm H}^1(\re)\right\},
    $$
    we find that the mini-max principle implies
    \bel{san13}
    {\rm Tr}\,{\bf 1}_{(-\infty, 2bq)}(H^{(q)}) \geq {\rm Tr}\,{\bf 1}_{(-\infty, 0)}(h_q).
    \ee
    Further, we could assume without loss of generality that the set $K$ in \eqref{san11} is bounded and axisymmetric, while the set $\tilde{K}$ in \eqref{san12} is radially symmetric.
    By the upper bounds in \eqref{san11}--\eqref{san12} and the mini-max principle, we have
    \bel{san14}
     {\rm Tr}\,{\bf 1}_{(-\infty, 0)}(h_q) \geq  {\rm Tr}\,{\bf 1}_{(-\infty, 0)}(\widetilde{h_q})
    \ee
    where $\widetilde{h_q} : = - \frac{d^2 }{d x_3^2} + \widetilde{v_q}$, and
    $$
    \widetilde{v_q}(x_3) : = -C\left\{
    \begin{array} {l}
     \int_0^{\infty} {\bf 1}_K(\rho,x_3) \phi_q(\rho)^2 \rho d\rho \quad \mbox{if \eqref{san11} holds true},\\
   \langle x_3\rangle^{-m_3} \int_0^{\infty} {\bf 1}_{\tilde{K}}(\rho) \phi_q(\rho)^2 \rho d\rho \quad \mbox{if \eqref{san12} holds true},
    \end{array}
    \right.
    \quad x_3 \in \re.
    $$
    Now note that the operator $\widetilde{h_q}$ has at least one negative eigenvalue   if \eqref{san11} holds true (see e.g. \cite[Lemma 5.2]{ahs}), and it has an infinite sequence of negative discrete eigenvalues accumulating at the origin if \eqref{san12} holds true (see e.g. \cite[Theorem XIII.82]{rs4}). Therefore, estimates \eqref{san15}--\eqref{san16} now follow from \eqref{san13}--\eqref{san14}.\\

{\sl Acknowledgments.}
G. Raikov thanks Fumio Hiroshima for the opportunity to give a talk on the results of the present article at the Conference Spectral and Scattering Theory and Related Topics, RIMS, Kyoto,
Japan, 14 - 17 December 2011.\\
J.-F. Bony and V. Bruneau
were partially supported by ANR-08-BLAN-0228.  G.
Raikov was partially supported by the Chilean Science Foundation
{\em Fondecyt} under Grant 1130591, and by {\em N\'ucleo Cient\'ifico
ICM} P07-027-F ``{\em Mathematical Theory of Quantum and Classical
Magnetic Systems''}.

{\sc Jean-Fran\c{c}ois Bony }\\
Universit\'e Bordeaux I, Institut de Math\'ematiques de
Bordeaux,\\
UMR CNRS 5251,  351, Cours de la Lib\'eration, 33405 Talence,
France\\
E-mail: bony@math.u-bordeaux1.fr\\

{\sc Vincent Bruneau }\\
Universit\'e Bordeaux I, Institut de Math\'ematiques de
Bordeaux,\\
UMR CNRS 5251,  351, Cours de la Lib\'eration, 33405 Talence,
France\\
E-mail: vbruneau@math.u-bordeaux1.fr\\

{\sc Georgi Raikov}\\
 Departamento de Matem\'atica, Facultad de
Matem\'aticas,\\ Pontificia Universidad Cat\'olica de Chile,
Vicu\~na Mackenna 4860, Santiago de Chile\\
E-mail: graikov@mat.puc.cl

\end{document}